\documentclass{llncs}
\usepackage{graphicx}   
\usepackage{amsmath} 
\usepackage{amssymb} 
\usepackage{url}

\begin{document}

\title{ Amari Functors and Dynamics in Gauge Structures}

\author{Michel Nguiffo Boyom\inst{1} \and Ahmed Zeglaoui\inst{1,2}} 

 \institute{IMAG\\ Institut Montpelli\'erain Alexander Grothendieck \\ Universit\'e de Montpellier. CC051, Place Eug\'ne Bataillon \\ F-34095 Montpellier, France\\ \email{boyom@math.univ-montp2.fr}\\
 {2} Laboratory of Algebra and Number Theory \\ Facult\'e de Math\'ematiques \\
 University of Science and Technology Houari Boumediene (USTHB) \\
 16111 Bab Ezzouar, Algeria \\
\email{ahmed.zeglaoui@gmail.com}}
\maketitle

\begin{abstract} 
We deal with finite dimensional differentiable manifolds. All items are concerned with are differentiable as well. The  class of differentiability is $C^\infty$. A metric structure in a vector bundle $E$ is a constant rank symmetric bilinear vector bundle homomorphism of $E\times E$ in the trivial bundle line bundle. We address the question whether a given gauge structure in $E$ is metric. That is the main concerns. We use generalized Amari functors of the information geometry for introducing two index functions defined in the moduli space of gauge structures in $E$. Beside we introduce a differential equation whose analysis allows to link the  new index functions just mentioned with the main concerns. We sketch applications in the differential geometry theory of statistics. \\

Reader interested in a former forum on the question whether a giving connection is metric are referred to appendix.

\end{abstract}
\section{Introduction}
A metric structure $(\mathbf{E},\mathbf{g})$ in a vector bundle $\mathbf{E}$ assigns to every fiber $\mathbf{E}_x$ a symmetric bilinear form $\mathbf{g}_x : \mathbf{E}_x \times \mathbf{E}_x \rightarrow \mathbb{R}$. Every finite rank vector bundle admits nondegenerate positive metric structures. One uses the paracompacity for constructing those positive regular metric structures. At another side every nondegenerate metric vector bundle $(\mathbf{E},\mathbf{g})$ admits metric gauge structures, viz gauge structures $(\mathbf{E},\mathbf{\nabla})$ subject to the requirement $\mathbf{\nabla}\mathbf{g} = 0$.
In a nondegenerate structure the values of the curvature tensor of a metric gauge structure $(\mathbf{E},\mathbf{\nabla})$ belong to the orthogonal sub-algebra $o(\mathbf{E}, \mathbf{\nabla})$ of the Lie algebra $\mathbb{G}(\mathbf{E})$. Arises the question whether a gauge structure $(\mathbf{E},\mathbf{\nabla})$ is a metric gauge structure in $\mathbf{E}$. \\
Our concern is to relate this existence question with some methods of the information geometry. In fact in the family $\nabla^\alpha$ of $\alpha$-connections in a non singular statistical model $[\mathbf{E},\pi,M,D,p]$ the 0-connection yields a metric gauge structure in $(TM,g)$. Here $g$ in the Fisher information of the statistical model as in \cite{ama:nag}, \cite{boyom1}. The question what about the cases $\alpha\neq 0$ deserves the attention. More generally arises the question when the pair $(\nabla, \nabla^{\star})$ in a statistical manifold $(M, g, \nabla, \nabla^{\star})$ is a pair of a metric gauge structures? Our aim is to address those questions in the general framework of finite rank real vector bundle over finite dimensional smooth manifolds. Our investigation involve two dynamics in the category $\mathfrak{Ga}(\mathbf{E})$. The first dynamic is the natural action of the gauge group $\mathcal{G}(\mathbf{E})$. The second is the action of the infinitely generated Coxeter group generated by the family $\mathfrak{Me}(\mathbf{E})$ of regular  metric structures $(\mathbf{E}, \mathbf{g})$. This second dynamic is derived from Amari functors.
\section{The gauge dynamic in $\mathfrak{Ga}(\mathbf{E})$}
\subsection{The gauge group of a vector bundle}
Let $\mathbf{E^{\star}}$ be the dual vector bundle of $\mathbf{E}$. Throughout this section $2$ we go to identify the vector bundles $\mathbf{E^{\star}}\otimes \mathbf{E}$ and $Hom(\mathbf{E}, \mathbf{E})$. Actually $Hom(\mathbf{E},\mathbf{E})$ is the vector bundle of vector bundle homomorphisms from $\mathbf{E}$ to $\mathbf{E}$. The sheaf of sections of $\mathbf{E^{\star}}\otimes \mathbf{E}$ is denoted by $\bf\mathbb{G}(\mathbf{E})$. This $\bf\mathbb{G}(\mathbf{E})$ is a Lie algebra sheaf bracket is defined by
$$
\left(\phi ,\psi \right) \longmapsto \left[
\phi ,\psi \right] = \phi \circ \psi -\psi \circ \phi.
$$
Actually $\mathbf{E^{\star}}\otimes \mathbf{E}$ is a Lie algebras bundle. It is called the Lie algebra of infinitesimals gauge transformations. The sheaf of inversible sections of $\mathbf{E}^\star\otimes\mathbf{E}$ is denoted by $\mathcal{G}(\mathbf{E})$. This $\mathcal{G}(\mathbf{E})$ is a Lie groups sheaf whose composition is the composition of applications of $\mathbf{E}$ in $\mathbf{E}$. Elements of $\mathcal{G}(\mathbf{E})$ are called gauge transformations of the vector bundle $\mathbf{E}$.
Consequently the set $\mathcal{G}_{x}(\mathbf{E}) \subset Hom(\mathbf{E}_{x}, \mathbf{E}_{x})$ is nothing but the Lie group $GL(\mathbf{E}_{x})$. This $\mathcal{G}(\mathbf{E})$ is the seheaf of sections of the Lie groups bundle $\tilde{\mathbf{E^{\star}}\otimes \mathbf{E}} \subset \mathbf{E^{\star}}\otimes \mathbf{E}$.
We abuse by calling $\mathcal{G}(\mathbf{E})$ and $\tilde{\mathbf{E^{\star}}\otimes \mathbf{E}}$ the gauge group of the vector bundle $\mathbf{E}$. 
\subsection{Gauge structures in a vector bundle $\mathbf{E}$}
A gauge structure in a vector bundle $\mathbf{E}$ is a pair $(\mathbf{E}, \mathbf{\nabla})$ where $\mathbf{\nabla}$ is a Koszul connection in $\mathbf{E}$. The set of gauge structures is denoted by $\mathfrak{Ga}(\mathbf{E})$.
We define the action of the gauge group in $\mathfrak{Ga}(\mathbf{E})$ as it follows
$$
\mathcal{G}(\mathbf{E}) \times \mathfrak{Ga}(\mathbf{E}) \longrightarrow \mathfrak{Ga}(\mathbf{E}),
$$
$$
\phi^{\star}(\mathbf{E}, \mathbf{\nabla}) = (\mathbf{E}, \phi^{\star}\mathbf{\nabla}).
$$
The Koszul connection $\phi^{\star}\mathbf{\nabla}$ is defined by
$$
(\phi^{\star}\mathbf{\nabla})_{X}s = \phi(\mathbf{\nabla}_{X}\phi^{-1}s)
$$
for all $s \in \mathfrak{G}a(\mathbf{E})$ and all vectors field $X$ on $M$.
We denoted the gauge moduli space by $Ga(\mathbf{E})$, viz
$$
Ga(\mathbf{E}) = \frac{\mathfrak{Ga}(\mathbf{E})}{\mathcal{G}(\mathbf{E})}
$$
\subsection{The equation $FE(\mathbf{\nabla}\mathbf{\nabla}^{\star})$}
Inspired by the appendix to \cite{boyom1} and by \cite{boyom3} and by  we define a map from pairs of gauge structures in the space of differential operators $DO(Eo^*\otimes E,T^* M\otimes E^*\otimes E).$
To every pair of gauge structures $\left[(\mathbf{E}, \mathbf{\nabla}) , (\mathbf{E}, \mathbf{\nabla}^{\star})\right]$ we introduce the first order differential operator $D^{\mathbf{\nabla}\mathbf{\nabla}^{\star}}$ of $\mathbf{E^{\star}}\otimes \mathbf{E}$ in $T^{\star}M \otimes\mathbf{E^{\star}}\otimes \mathbf{E}$ which is defined as it follows 
$$
D^{\nabla\nabla^*}(\phi)(X,s) = \nabla^*_X(\phi(s)) - \phi(\nabla_Xs)
$$ 
for all $s$ and for all vector fields $X$.
Assume the rank of $\mathbf{E}$\ is equal to $r$ and the dimension of $M$ is
equal to $m$. Assume $\left( x_{i}\right) $ is a system of local coordinate
functions defined in an open subset $U \subset M$ and $\left( s_{\alpha }\right)$ is
a basis of local sections of $\mathbf{E}$ defined in $U$. We set
\[
\begin{array}{ccccc}
\mathbf{\nabla }_{\partial _{x_i}}s_{\alpha }=\sum_{\beta }\Gamma _{i:\alpha
}^{\beta }s_{\beta }, &  & \mathbf{\nabla }_{\partial _{x_i}}^{\star
}s_{\alpha }=\sum_{\beta }\Gamma _{i:\alpha }^{\star \beta }s_{\beta } &
\text{and} & \phi \left( s_{\alpha }\right) =\sum_{\beta }\phi _{\alpha
}^{\beta }s_{\beta }.
\end{array}
\]
Our concern is the analysis of system of partial derivative equations%
\[
\left[ FE\left( \mathbf{\nabla \nabla }^{\star }\right) \right] _{i:\alpha
}^{\gamma }:\frac{\partial \phi _{\alpha }^{\beta }}{\partial x_{i}}%
+\sum_{\beta =1}^{r}\left\{ \phi _{\alpha }^{\beta }\Gamma _{i:\beta
}^{\star \gamma }-\phi _{\beta }^{\gamma }\Gamma _{i:\alpha }^{\beta
}\right\} =0.
\]

When we deal with the vector tangent bundles the differential operator $D^{\nabla\nabla^*}$ plays many outstanding roles in the global analysis of the base manifold \cite{boyom3}. In general though every vector bundle admits positive metric structures this  same claim is far from being true for symplectic structure and for positive signature metric structures. We aim at linking those open problems with the differential equation $FE(\nabla\nabla^*)$. \\
The sheaf of germs of solutions to $FE(\mathbf{\nabla} \mathbf{\nabla}^{\star})$ is denoted by $\mathcal{J}_{\mathbf{\nabla}\mathbf{\nabla}^{\star}}(\mathbf{E})$.
\section{The metric dynamics in $\mathfrak{Ga}(\mathbf{E})$}
\subsection{The Amari functors in the category $\mathfrak{Ga}(\mathbf{E})$}
Without the express statement of the contrary a metric structure in a vector bundle $\mathbf{E}$ is a constant rank symmetric bilinear vector bundle homomorphism $\mathbf{g}$ of $\mathbf{E} \otimes \mathbf{E}$ in $\tilde{\mathbb{R}}$. Such a metric structure is denoted by $(\mathbf{E}, \mathbf{g})$. A nondegenerate metric structure is called regular, otherwise it is called singular. The category of regular metric structures in $\mathbf{E}$ is denoted by $\mathfrak{Me}(\mathbf{E})$. \\
Henceforth the concern is the dynamic
\[
\begin{array}{ccc}
\mathcal{G}(\mathbf{E}) \times \mathfrak{Me}(\mathbf{E}) &
\longrightarrow  & \mathfrak{Me}\mathcal{(}\mathbf{E}\mathcal{)} \\
\left( \phi ,\left( \mathbf{E},\mathbf{g}\right) \right)  & \longmapsto
& \left( E,\phi _{\star } \mathbf{g} \right).
\end{array}
\]
Here the metric $\phi _{\star }\mathbf{g} $ is defined by
\[
\phi _{\star }\mathbf{g} (s,s^{\prime }) = \mathbf{g(}\phi
^{-1}(s),\phi ^{-1}(s^{\prime })\mathbf{).}
\]
This leads to the moduli space of regular metric structures in a vector bundle $\mathbf{E}$
$$
Me(\mathbf{E}) = \frac{\mathfrak{Me}(\mathbf{E})}{\mathcal{G}(\mathbf{E})}.
$$
A gauge structure $(\mathbf{E}, \mathbf{\nabla})$ is called metric if there exist a metric structure $(\mathbf{E}, \mathbf{g})$ subject to the requirement $\mathbf{\nabla}\mathbf{g} = 0$.  \\
We consider the functor $\mathfrak{Me}(\mathbf{E}) \times \mathfrak{Ga}(\mathbf{E}) \rightarrow \mathfrak{Ga}(\mathbf{E})$ which is defined by 
$$
\left[ (\mathbf{E}, \mathbf{g}) , (\mathbf{E}, \mathbf{\nabla}) \longmapsto (\mathbf{E}, \mathbf{g.\nabla}) \right].$$ Here the Koszul connection $\mathbf{g.\nabla}$
is defined by
$$
\mathbf{g}(\mathbf{g.\nabla}_{X}s, s^{\prime}) = X(\mathbf{g}(s, s^{\prime}))-\mathbf{g}(s, \mathbf{\nabla}_{X}s^{\prime}).
$$
The functor just mentioned is called the general Amari functor of the vector bundle $\mathbf{E}$. According to \cite{boyom3}, the general Amari functor yield two restrictions :
\begin{equation}
\begin{array}{ccc}
\left\{ \mathbf{g}\right\} \times \mathfrak{Ga}\left( \mathbf{E},\right)  &
\longrightarrow  & \mathfrak{Ga}\left( \mathbf{E}\right)  \\
\mathbf{\nabla } & \longmapsto  & \mathbf{g.\nabla }
\end{array}
\end{equation}
\begin{equation}
\begin{array}{ccc}
\mathfrak{Me}\left( \mathbf{E}\right) \times \left\{ \mathbf{\nabla }\right\}
& \longrightarrow  & \mathfrak{Ga}\left( \mathbf{E}\right)  \\
\mathbf{g} & \longmapsto  & \mathbf{g.\nabla }
\end{array}
\end{equation}
The restriction $(1)$ is called the metric Amari functor of the gauge structure $(\mathbf{E},\mathbf{\nabla})$. The restriction $(2)$ is called the gauge Amari functor of the metric vector bundle $(\mathbf{E},\mathbf{g}) $ . \\
We observe that $\mathbf{\nabla}\mathbf{g} = 0$ if and only if $\mathbf{g.\nabla} =  \mathbf{\nabla}$. 
The restriction $(1)$ gives rise to the involution of $\mathfrak{Ga}( \mathbf{E})$ : $\mathbf{\nabla} \rightarrow \mathbf{g}.\mathbf{\nabla}$. In other words $\mathbf{g}.(\mathbf{g.\nabla}) = \mathbf{\nabla}$ for all $(\mathbf{E}, \mathbf{\nabla}) \in \mathfrak{Ga}(\mathbf{E})$. 
In general the question whether an involution admits fixed points has negative answers. In the framework we are concerned with  every involution defined by a regular metric structure has fixed points formed by metric gauge structures in $(\mathbf{E}, \mathbf{g})$.\\
The dynamics
\[
\begin{array}{ccc}
\mathcal{G}\left( \mathbf{E}\right) \times \mathfrak{Ga}\left( \mathbf{E}%
\right)  & \longrightarrow  & \mathfrak{Ga}\left( \mathbf{E}\right)  \\
\left( \phi ,\mathbf{\nabla }\right)  & \longmapsto  & \phi ^{\star
} \mathbf{\nabla }
\end{array}%
\]%
\[
\begin{array}{ccc}
\mathcal{G}\left( \mathbf{E}\right) \times \mathfrak{Me}\left( \mathbf{E}%
\right)  & \longrightarrow  & \mathfrak{Me}\left( \mathbf{E}\right)  \\
\left( \phi ,\mathbf{g}\right)  & \longmapsto  & \phi _{\star }
\mathbf{g}
\end{array}%
\]%
are linked with the metric Amari functor by the formula
\[
\phi ^{\star } \mathbf{g.\nabla } = \phi _{\star }
\mathbf{g} . \phi ^{\star } \mathbf{\nabla }.
\]
We go to introduce the metric dynamics in $\mathfrak{Ga}(\mathbf{E})$. The abstract group of all isomorphisms of $\mathfrak{Ga}(\mathbf{E})$ is denoted by $ISO(\mathfrak{Ga}(\mathbf{E}))$. By the metric Amari functor every regular metric structure $(\mathbf{E}, \mathbf{g})$ yields the involution $(\mathbf{E}, \mathbf{\nabla}) \rightarrow (\mathbf{E}, \mathbf{g.\nabla})$.
The subgroup of $ISO(\mathfrak{Ga}(\mathbf{E}))$ which is generated by all regular metric structures in $\mathbf{E}$ is denoted by $\mathcal{G}m(\mathbf{E})$. This group $\mathcal{G}m(\mathbf{E})$ looks like an infinitely generated Coxeter group.
Using this analogy we call $\mathcal{G}m(\mathbf{E})$ the metric Coxeter group of $\mathfrak{Ga}(\mathbf{E})$.
For instance every metric structure $(\mathbf{E}, \mathbf{g})$ generates a dihedral group of order $2$.
\subsection{The quasi-commutativity property of the metric dynamic and the gauge dynamic}
At the present step we are dealing with both the gauge dynamic
\begin{equation}
\begin{array}{ccc}
\mathcal{G}\left( \mathbf{E}\right) \times \mathfrak{Ga}\left( \mathbf{E}
\right)  & \longrightarrow  & \mathfrak{Ga}\left( \mathbf{E}\right)  \\
\left( \phi ,\mathbf{\nabla }\right)  & \longmapsto  & \phi ^{\star
} \mathbf{\nabla }
\end{array}
\end{equation}
and the metric dynamic
\begin{equation}
\begin{array}{ccc}
\mathcal{G}m\left( \mathbf{E}\right) \times \mathfrak{Ga}\left( \mathbf{E}
\right)  & \longrightarrow  & \mathfrak{Ga}\left( \mathbf{E}\right)  \\
\left( \gamma, \mathbf{\nabla }\right)  & \longmapsto  & \gamma .\mathbf{
\nabla }
\end{array}
\end{equation}
What we call the quasi commutativity property of $(1)$ and $(2)$ is the link
\[
\phi ^{\star } \mathbf{g.\nabla} = \phi _{\star }
\mathbf{g} .\phi ^{\star } \mathbf{\nabla }
\].
We consider two regular metric structures $(\mathbf{E}, \mathbf{g}^{0})$ and $(\mathbf{E}, \mathbf{g})$, There exists a unique $\phi\in \mathcal{G}(\mathbf{E})$ subject to the requirement
$$
\mathbf{g}^{0}(s, s^{\prime}) = \mathbf{g}(\phi(s), s^{\prime}).
$$
By direct calculations one sees that for every gauge structure $(\mathbf{E}, \mathbf{\nabla}) $ one has
$$
\mathbf{g}.\mathbf{\nabla} = \phi^{\star}(\mathbf{g}^{0}.\mathbf{\nabla}).
$$
The quasi-commutativity property shows that every regular metric structure acts in the moduli space $Ga(\mathbf{E})$. Further the gauge orbit $ \mathcal{G}(\mathbf{E})(\mathbf{g.\nabla})$ does not depend on the choice of the regular metric structure $(\mathbf{E},\mathbf{g})$. Thus the metric Coxeter group $\mathcal{G}m(\mathbf{E})$ acts in the moduli space $Ga(\mathbf{E}).$  When there is no risk of confusion the orbit of $[\mathbf{\nabla}] \in Ga(\mathbf{E})$ is denoted by 
$\mathcal{G}m.[\mathbf{\nabla}]$ while its stabilizer subgroup is denoted by $\mathcal{G}m_{[\mathbf{\nabla}]}$. Consequently one has
\begin{proposition} 
The index of every stabilizer subgroup $\mathcal{G}m_{[\mathbf{\nabla}]} \subset \mathcal{G}m(\mathbf{E})$ is equal to $1$ or to $2$.
\end{proposition}
We go to rephrase Proposition 1 versus the orbits of the metric Coxeter group in the moduli space $ Ga(\mathbf{E})$.
\begin{proposition}
For every orbit $\mathcal{G}m(\mathbf{E}).[\mathbf{\nabla}]$  cardinal $\sharp \left( \mathcal{G}m(\mathbf{E}).[\mathbf{\nabla}] \right) \in \left\{ 1, 2 \right\} $
\end{proposition}
\subsection{The metric index function.}
The  concern is the metric dynamic 
\begin{equation}
\begin{array}{ccc}
\mathcal{G}m\left( \mathbf{E}\right) \times \mathfrak{Ga}\left( \mathbf{E}
\right)  & \longrightarrow  & \mathfrak{Ga}\left( \mathbf{E}\right)  \\
\left( \gamma \mathbf{,\nabla }\right)  & \longmapsto  & \gamma .\mathbf{
\nabla }
\end{array}
\end{equation}
The length of $\gamma \in \mathcal{G}m(\mathbf{E} )$ is denoted by $\mathfrak{l}(\gamma).$ It is defined as it follows
$$
\mathfrak{l}(\gamma) = \min \left \{ p \in \mathbb{N} : \gamma = \mathbf{g}_{1}\mathbf{g}_{2}\ldots \mathbf{g}_{p},\quad \mathbf{g}_{j}\in \mathfrak{M}e(\mathbf{E}) \right \}
$$
For every gauge structure $(\mathbf{E}, \mathbf{\nabla})$ the metric index of $(\mathbf{E},\mathbf{\nabla})$ is defined by
$$ind(\mathbf{\nabla}) = \min_{{\gamma \in \mathcal{G}^*m_{\mathbf{\nabla}}}} \left \{ \mathfrak{l}(\gamma) - 1\right\}.
$$
Here $\mathcal{G}^*m_{\mathbf{\nabla}}$ stands for the subset formed of elements of the isotropy subgroup that differ from the unit element. The flowing statement is a  straightforward consequence of the quasi-commutativity property.  
\begin{lemma}The non negative integer $ind(\mathbf{\nabla})$ is a gauge invariant.
\end{lemma}
Consequently we go to encode every orbit $ [\mathbf{\nabla}] = \mathcal{G}(\mathbf{E})^{\star}\mathbf{\nabla}$ with metric index  $ind([\mathbf{\nabla}]) =  ind(\mathbf{\nabla})$.
\begin{definition}
By Lemma 1 we get the metric index function
$$\mathcal{G}a({\mathbf{E}}) \ni [\mathbf{\nabla}]\rightarrow ind([\mathbf{\nabla}]) \in \mathbb{Z}$$
\end{definition}
\subsection{The gauge index function.}
We consider the general Amari functor
$$
\begin{array}{ccc}
\mathfrak{Me}\left( \mathbf{E}\right) \times \mathfrak{Ga}\left( \mathbf{E}
\right)  & \longrightarrow  & \mathfrak{Ga}\left( \mathbf{E}\right)  \\
\left( \mathbf{g,\nabla }\right)  & \longmapsto  & \mathbf{g.\nabla }
\end{array}
$$
For convenience we set $\mathbf{\nabla}^{\mathbf{g}} = \mathbf{g}.\mathbf{\nabla}$.
Therefore to a pair $[(\mathbf{E}, \mathbf{g}), (\mathbf{E}, \mathbf{\nabla})]$ we assign the differential equation $FE(\mathbf{\nabla}\mathbf{\nabla}^{\mathbf{g}})$. The sheaf of solutions to $FE(\mathbf{\nabla}\mathbf{\nabla}^{\mathbf{g}})$ is denoted by $\mathcal{J}_{\mathbf{\nabla}\mathbf{\nabla}^{\mathbf{g}}}(\mathbf{E})$.
We go to perform a formalism which is developed in \cite{boyom2}. See also \cite{boyom3} for the case of tangent bundles of a manifolds.
The concerns are metric structures in vector bundles. We recall that a singular metric structure in $\mathbf{E}$ is a constant rank degenerate symmetric bilinear vector bundles homomorphism $\mathbf{g} : \mathbf{E} \times \mathbf{E} \rightarrow \tilde{\mathbb{R}}$.
Let $(\mathbf{E},\mathbf{g})$ be a regular metric structure. We pose $\mathbf{\nabla}^{\mathbf{g}} = \mathbf{g}.\mathbf{\nabla}$.
\newline
For every $\phi \in \mathbb{G}(\mathbf{E})$ there exists a unique pair $(\Phi,\Phi^* ) \subset \mathbb{G}(\mathbf{E})$ subject to the the following requirements
\begin{equation}
\mathbf{g}(\Phi(s), s^{\prime}) = \frac{1}{2}\left[ \mathbf{g}(\phi(s), s^{\prime})+
\mathbf{g}(s, \phi(s^{\prime}))
 \right],
\end{equation}
\begin{equation}
\mathbf{g}(\Phi^{\star}(s), s^{\prime}) = \frac{1}{2}\left[ \mathbf{g}(\phi(s), s^{\prime})-
\mathbf{g}(s, \phi(s^{\prime}))
 \right]
\end{equation}
We put $q(s, s^{\prime}) = \mathbf{g}(\Phi(s), s^{\prime})$ and $\omega (s, s^{\prime}) = \mathbf{g}(\Phi^\star(s), s^{\prime})$.
\begin{proposition} (\cite{boyom2}) 
If $\phi$ is a solution to $FE(\mathbf{\nabla}\mathbf{\nabla}^{\star})$ then $\Phi$ and $\Phi^{\star}$ are solutions to $FE(\mathbf{\nabla}\mathbf{\nabla}^{\star})$. Furthermore,
$$\mathbf{\nabla}q = 0,$$
$$\mathbf{\nabla}\omega = 0. $$
\end{proposition}
By the virtue of the proposition $3$ one has $rank(\Phi) = Constant$  and $rank(\Phi^{\star}) = Constant$. 
\begin{corollary}
We assume that the  regular metric structure $(\mathbf{E}, \mathbf{g})$ is positive definite then we have 
\begin{equation}
\mathbf{E} = \ker(\Phi) \oplus Im(\Phi)
\end{equation}
\begin{equation}
\mathbf{E} = \ker(\Phi^{\star}) \oplus Im(\Phi^{\star})
\end{equation}
Further one has the following gauge reductions
\begin{equation}
(\ker(\Phi), \mathbf{\nabla}) \subset (\mathbf{E}, \mathbf{\nabla}),
\end{equation}
\begin{equation}
(Im(\Phi), \mathbf{\nabla}^{\mathbf{g}}) \subset (\mathbf{E}, \mathbf{\nabla}^{\mathbf{g}}),
\end{equation}
\begin{equation}
(\ker(\Phi^{\star}), \mathbf{\nabla}) \subset (\mathbf{E}, \mathbf{\nabla}),
\end{equation}
\begin{equation}
(Im(\Phi^{\star}), \mathbf{\nabla}^{\mathbf{g}}) \subset (\mathbf{E}, \mathbf{\nabla}^{\mathbf{g}}).
\end{equation}
\end{corollary}
\begin{corollary} Assume that $(\mathbf{E}, \mathbf{g}, \mathbf{\nabla}, \mathbf{\nabla}^{\mathbf{g}})$ is the vector bundle versus of a statistical manifold $(M, \mathbf{g}, \mathbf{\nabla}, \mathbf{\nabla}^{\mathbf{g}})$ here $\mathbf{\nabla} = $. Then $(10, 11, 12, 13)$  is (quasi) $4$-web in the base  manifold $M$.
\end{corollary}
Given a metric vector bundle $(\mathbf{E}, \mathbf{g})$ and gauge structure $(\mathbf{E}, \mathbf{\nabla})$. The triple $(\mathbf{E}, \mathbf{g}, \mathbf{\nabla})$ is called special if the differential equation $FE(\mathbf{\nabla}\mathbf{\nabla}^{\mathbf{g}})$ has non trivial solutions.
We deduce from Corollary $2$  that every special statistical manifold supports a canonical (quasi) $4$-web, viz $4$ foliations in (quasi) general position. \\
Before pursing we remark that among formalisms introduce in \cite{boyom3}, many (of them) walk in the category of vector bundles. We go to perform this remark. To every special triple $(\mathbf{E}, \mathbf{g}, \mathbf{\nabla})$ we assign the function
\[
\begin{array}{ccc}
\mathcal{J}_{\mathbf{\nabla \nabla }^{\mathbf{g}}}\left( \mathbf{E}\right)
& \longrightarrow  & \mathcal{\mathbb{Z}} \\
\phi  & \longmapsto  & rank(\Phi )
\end{array}
\]
\textbf{Reminder} : The map $\Phi$ is the solution to $FE(\nabla\nabla^g)$ given by 
$\mathbf{g}(\Phi(s), s^{\prime}) = \frac{1}{2}\left[ \mathbf{g}(\phi(s), s^{\prime})+
\mathbf{g}(s, \phi(s^{\prime})) \right]$.
We define the following non negatives integers
\begin{equation}
s^{b}\left( \mathbf{\nabla },\mathbf{g}\right) = \min_{\phi \in \mathcal{J}%
_{\mathbf{\nabla \nabla }^{\mathbf{g}}}\left( \mathbf{E}\right)
}corank\left( \Phi \right) ,
\end{equation}
\begin{equation}
s^{b}\left( \mathbf{\nabla }\right) = \min_{\left( \mathbf{E,g}\right) \in
\mathfrak{Me}\left( \mathbf{E}\right) }s^{b}\left( \mathbf{\nabla },\mathbf{g}%
\right) .
\end{equation}
\begin{proposition}
The non negative integer $s^b(\nabla)$ is a gauge invariant   $\mathfrak{G}a(\mathbf{E})$, viz $s^b(\mathbf{\nabla}) = s^b(\Phi^\star \mathbf{\nabla})$ for all gauge transformation $\Phi$.
\end{proposition}
\begin{definition}
By Proposition $4$ we get the gauge index function
$$
\mathcal{G}a(\mathbf{E}) \ni [\mathbf{\nabla}] \rightarrow s^b([\mathbf{\nabla}]) \in \mathbb{Z} $$
\end{definition}
\section{The topological nature of the index functions}
\subsection{Index functions as characteristic obstruction}
According to \cite{aff:che}, every positive Riemannian foliation (nice singular metric in the tangent bundle of a smooth manifold) admits a unique symmetric metric connection. A combination of \cite{aff:che} and \cite{boyom1} shows that all those metrics are constructed using methods of the information geometry as in \cite{boyom1} (see the exact sequence $(16)$ below).
Remind that we are concerned with the question whether a gauge structure $(\mathbf{E}, \mathbf{\nabla})$ is metric.
By the virtue of \cite{aff:che} and \cite{boyom2} one has 
\begin{theorem} In a finite rank vector bundle $\mathbf{E}$ a gauge structure $(\mathbf{E},\mathbf{\nabla})$ is metric if and only if for some regular metric structure $(\mathbf{E}, \mathbf{g})$ the differential equation 
$FE(\nabla\nabla^\mathbf{g})$ admits non trivial solutions.
\end{theorem}
\begin{remark}   
If for some metric structure $(\mathbf{E},\mathbf{g}^0)$ the differential equation $FE(\nabla\nabla^{g^0})$ admits non trivial solutions then for every regular metric structure $(\mathbf{E}, \mathbf{g})$ the differential equation  $FE(\nabla\nabla^{g})$ admits non trivial solutions.\\
Hint : use the following the short exact sequence as in \cite{boyom2}
\begin{equation}
0 \longrightarrow \Omega_{2}^{\mathbf{\nabla}}(TM) \longrightarrow \mathcal{J}_{\mathbf{\nabla}\mathbf{\nabla^{g}}}(TM) \longrightarrow S_{2}^{\mathbf{\nabla}}(TM) \longrightarrow 0.
\end{equation}
\end{remark}
We recall that the concern is the question whether a gauge structure $(\mathbf{E}, \mathbf{\nabla})$ is metric. By the remark raised above, this question is linked with the solvability of differential equations $FE(\mathbf{\nabla}\mathbf{\nabla}^{\mathbf{g}})$ which locally is a system of linear PDE with non constant coefficients. Theorem 1 highlights the links of its solvability with the theory of Riemannian foliations which are objects of the differential topology. The key of those links are items of the information geometry. So giving $(\mathbf{E}, \mathbf{\nabla})$, the property of $(\mathbf{E}, \mathbf{\nabla})$ to be metric is equivalent to the property of $FE(\mathbf{\nabla}\mathbf{\nabla}^{\mathbf{g}})$ to admit non trivial solutions.
Henceforth, our aim is to relate the question just mentioned and the invariants $ind(\mathbf{\nabla})$ and $s^{b}(\mathbf{\nabla})$. We assume that $(\mathbf{E}, \mathbf{\nabla})$ is special.
\begin{theorem} In a gauge structure $(\mathbf{E}, \mathbf{\nabla})$, the following assertions are equivalent
\begin{enumerate}
\item The gauge structure $(\mathbf{E}, \mathbf{\nabla})$ is regularly special.
\item The metric index function vanishes at $[\mathbf{\nabla}] \in \mathcal{G}a(\mathbf{E})$ i.e. $ind([\mathbf{\nabla}]) = 0.$
\item The gauge index function vanishes at $[\mathbf{\nabla}] \in \mathcal{G}a(\mathbf{E})$ i.e. $s^{b}([\mathbf{\nabla}]) = 0$.
\end{enumerate}
\end{theorem}
\begin{theorem}
A gauge structure $(\mathbf{E},\mathbf{\nabla})$ is regularly metric if and only if $(\mathbf{E},\mathbf{g.\nabla})$ is regularly metric for all regular metric structure $(\mathbf{E},\mathbf{g})$
\end{theorem}
By theorem $1$ both $ind(\mathbf{\nabla})$ and $s^{b}(\mathbf{\nabla})$ are characteristic obstructions to $(\mathbf{E}, \mathbf{\nabla})$ being regularly special.
We have no relevant interpretation of the case $ind(\mathbf{\nabla}) \neq 0$.
Regarding the case $s^{b}(\mathbf{\nabla}) \neq 0$, we have
\begin{proposition}
Let $(\mathbf{E}, \mathbf{\nabla})$ be a gauge structure with $s^{b}(\mathbf{\nabla}) \neq 0$. Then there exists a metric structure $(\mathbf{E}, \mathbf{g})$ such subject to the following requirement : $rank(\mathbf{g}) = s^{b}(\mathbf{\nabla})$, further $\mathbf{g}$ is optimal for those requirement, viz every $\mathbf{\nabla}$-parallel metric structure $(\mathbf{E}, \mathbf{g})$ has rank smaller than $s^b(\mathbf{\nabla})$
\end{proposition}
\subsection{Applications to the statistical geometry}
\begin{theorem}
Let $\left\{\mathbf{\nabla}^\alpha \right\}$ be the family of 
$\alpha$-connections of a statistical manifold. If $\nabla^\alpha$ is regularly metric for all of the positive real numbers $\alpha$ then all of the $\alpha$-connections are regularly metric.
\end{theorem}

\textbf{Appendix : When can a Connection Induce a Riemannian Metric for which it is the Levi-Civita Connection?}

https://mathoverflow.net/questions/54434/when-can-a-connection-induce-a-riemannian-metric-for-which-it-is-the-levi-civita

\end{document}